\documentclass{amsart}
\usepackage{verbatim,amssymb,amsmath,amscd,latexsym,amsbsy,mathrsfs} \input{xy}
\xyoption{all}

\newtheorem{thm}{Theorem}[section] \newtheorem{pro}[thm]{Proposition}
 \newtheorem{lm}[thm]{Lemma}
\newtheorem{cor}[thm]{Corollary}

\theoremstyle{remark}

\theoremstyle{definition}

\DeclareMathAlphabet{\mathpzc}{OT1}{pzc}{m}{it}
\DeclareMathOperator*{\spec}{Spec} \DeclareMathOperator*{\Hom}{Hom}
\DeclareMathOperator*{\card}{card} \DeclareMathOperator*{\im}{Im}
\DeclareMathOperator*{\Aut}{Aut} \DeclareMathOperator*{\QF}{frac}
\DeclareMathOperator*{\Gal}{Gal} \DeclareMathOperator*{\Ind}{Ind}
\DeclareMathOperator*{\cd}{cd}

\DeclareMathOperator*{\image}{im}
\DeclareMathOperator*{\into}{\hookrightarrow}
\newcommand{\chash}{\mathcal{\#}} 
 \newcommand{\QQ}{\mathbb{Q}}
\newcommand{\ZZ}{\mathbb{Z}} \newcommand{\Aff}{\mathbb{A}}
\newcommand{\PP}{\mathbb{P}} 
\newcommand{\SSS}{\mathbb{S}} 
\newcommand{\cB}{\mathcal{B}} \newcommand{\cF}{\mathcal{F}}
\newcommand{\cAP}{\mathcal{AM}} \newcommand{\cC}{\mathcal{C}}
\newcommand{\cA}{\mathcal{A}} \newcommand{\cSP}{\mathcal{SM}}
\newcommand{\cO}{\mathcal{O}} \newcommand{\cP}{\mathcal{M}}
\newcommand{\cp}{\mathpzc{p}}
 
\newcommand{\cK}{\mathcal{K}} 
\newcommand{\cG}{\mathcal{G}} 
\newcommand{\dirlim} {\displaystyle \lim_{\longrightarrow}}
\newcommand{\invlim} {\displaystyle \lim_{\longleftarrow}}
\newcommand{\D}{\displaystyle}

\begin{document}
\title{Fundamental Group in nonzero characteristic} \thanks{
  Department of Mathematics, Purdue University.  }  \author{Manish
  Kumar}
\begin{abstract}
  A proof of freeness of the commutator subgroup of the fundamental
  group of a smooth irreducible affine curve over a countable algebraically 
  closed field of nonzero characteristic. A description of the 
  abelianizations of the fundamental groups of affine curves over 
  an algebraically closed field of nonzero characteristic is also given.
\end{abstract}
\date{} \maketitle

\section{Introduction}
The algebraic fundamental group of smooth curves over an algebraically
closed field of characteristic zero is a well understood object,
thanks to Grothendieck's Riemann existence theorem [SGAI, XIII,
Corollary 2.12, page 392] . But if the characteristic of the base
field is $p>0$ and the curve is affine then there may be wild
ramification so computing the algebraic fundamental group is not as
simple. Though Grothendieck's theorem gives a description of the
prime-to-$p$ part of the fundamental group, which is analogous to the
characteristic zero case.  But the structure of the whole group is
still elusive in spite of the fact that all the finite quotients of
this group are now known. A necessary and sufficient condition for a
group
to be a finite quotient of the fundamental group was conjectured by
Abhyankar (see the theorem below) and was proved by Raynaud [Ray]
(in the case of the affine line) and Harbater [Ha1] (for arbitrary
smooth affine curves). For a finite group $G$ and a prime number $p$,
let $p(G)$ denote the subgroup of $G$ generated by all the $p$-Sylow
subgroups. $p(G)$ is called the \emph{quasi-$p$} part of $G$.

\begin{thm}({\bf Raynaud, Harbater})
  Let $C$ be a smooth projective curve of genus $g$ over an
  algebraically closed field of characteristic $p>0$ and for some
  $n\ge 0$, let $x_0,..,x_n$ be some points on $C$. Then a finite
  group $G$ is a quotient of the fundamental group $\pi_1(C\setminus
  \{x_0,..,x_n\})$ if and only if $G/p(G)$ is generated by $2g+n$
  elements. In particular a finite group $G$ is a quotient of
  $\pi_1(\Aff ^1)$ if and only if $G=p(G)$, i.e., $G$ is a quasi-$p$
  group.
\end{thm}

The ``if part'' of the above theorem is the nontrivial part, the
``only if part'' was proved long back by Grothendieck.

From now on we shall assume that the characteristic of the base field
is $p>0$.  Consider the following exact sequence for the
fundamental group of a smooth affine curve $C$.
$$1\rightarrow \pi_1 ^c (C) \rightarrow \pi_1(C) \rightarrow \pi_1
^{ab}(C) \rightarrow 1$$
where $\pi_1 ^c(C)$ and $\pi_1 ^{ab}(C)$ are
the commutator subgroup and the abelianization of the fundamental
group $\pi_1(C)$ of C, respectively. In this paper we give a
description of the abelianization and show that the commutator
subgroup is a free profinite group. 
The result on the commutator subgroup falls
into the league of the so called Shafarevich conjecture for global
fields. Recall that the Shafarevich conjecture says that the
commutator subgroup of the absolute Galois group of the rationals
$\QQ$ is free. David Harbater ([Ha7]), Florian Pop ([Pop]) and later 
Dan Haran and Moshe Jarden ([HJ]) have shown, using
different patching methods, that the absolute Galois group of the function
field of a curve over an algebraically closed field is free. More Shafarevich
conjecture type results have been proved in [Ha8].

The second section of this paper covers definitions and notations. In
the third section we give a description of the $p$-part of the
abelianization of the algebraic fundamental group of any normal affine
algebraic variety over an algebraically closed field in terms of Witt
vectors. We deduce the fact that the abelianization of the algebraic
fundamental group determines $W_n(A)/P(W_n(A))$ as a group (see
Corollary 3.6) where $W_n(A)$ is the ring of finite Witt vectors over
the coordinate ring $A$ of the affine curve under consideration and
$P$ is the additive group endomorphism of $W_n(A)$ which sends
$(a_1,..,a_n)$ to $(a_1^p,..,a_n^p) - (a_1,..,a_n)$ (here $''-''$ is
subtraction in the Witt ring).  It is conjectured by Harbater that the
algebraic fundamental group should determine $A$ as a ring.

The rest of the paper is devoted to proving that the commutator
subgroup of the fundamental group of a smooth irreducible affine curve over a countable
algebraically closed field $k$ is a free profinite group of countable
rank. The fourth section consists of some group
theory results and a result on projectivity of the commutator
subgroup. These results allow us to reduce the problem to finding
proper solution for any split embedding problem with perfect quasi-$p$
group as the kernel, abelian $p$-group as kernel and prime-to-$p$
group as kernel (for definitions see Section $1$).

The fifth section is on finding solutions for all quasi-$p$ perfect
embedding problems and abelian $p$-group embedding problems. This is
relatively simple. 

The sixth section is the longest section and is devoted to finding
proper solutions for prime-to-$p$ embedding problems. Methods
on formal patching developed by Harbater (see [Ha1] and [Ha2]) and
their mild generalizations have been used in this section to prove the
desired result.

I would like to thank my advisor, Prof. Donu Arapura, for his guidance
and numerous suggestions in developing the theory and verification of
the proofs. I am also grateful to Prof David Harbater for his useful 
suggestions and comments which helped me obtain the result on the
commutator subgroup in this generality, i.e., for any smooth 
affine irreducible curve.

\section{Definitions and notations}

An embedding problem consists of surjections, $\phi:\pi \rightarrow G$
and
$\alpha:\Gamma \rightarrow G$         \\
$$
\xymatrix{
  &          &                        &\pi \ar @{-->}[dl]_{\psi} \ar[d]^{\phi}\\
  1\ar[r] & H \ar[r] & \Gamma \ar[r]_{\alpha} & G \ar[r] \ar[d] & 1\\
  & & & 1 }
$$
where $G$, $H$, $\Gamma$ and $\pi$ are groups and $H=\ker(\alpha)$.
It is said to have a \emph{weak solution} if there exists a group
homomorphism $\psi$ which makes the diagram commutative, i.e., $\alpha
\circ \psi =\phi$. Moreover, if $\psi$ is an epimorphism then it is
said to have a \emph{proper} solution. It is said to be a finite
embedding problem if $\Gamma$ is finite. All the embedding problems
considered here will be assumed to be finite.  It is said to be a
split embedding problem if there exists a group homomorphism from $G$
to $\Gamma$ which is a right inverse of $\alpha$. It is called a
quasi-$p$ embedding problem if $H$ is a quasi-$p$ group, i.e., $H$ is
generated by its Sylow $p$-subgroups and similarly it is called a
prime-to-$p$ embedding problem if $H$ is a prime-to-$p$ group, i.e.,
order of $H$ is prime to $p$. $H$ will sometimes be referred to as the kernel
of the embedding problem.

A profinite group is called \emph{free} if it is a profinite
completion of a free group. A \emph{generating set} of a profinite
group $\pi$ is a subset $I$ so that the closure of the group generated
by $I$ is the whole group $\pi$. The \emph{rank} of a profinite group
is the minimum of the cardinality of a generating set.

For a ring $R$, let $\QF(R)$ denote the total ring of quotients of
$R$. A ring extension $R \subset S$ is said to be \emph{generically
  separable} if $R$ is a domain, $\QF(S)$ is separable extension of
$\QF(R)$ and no nonzero element of $R$ becomes a zero divisor in $S$.
 
As in [Ha2], a morphism of schemes, $\Phi: Y\rightarrow X$, is said to
be a \emph{cover} if $\Phi$ is finite and generically separable, i.e.,
$X$ can be covered by affine open subset $U=\spec(R)$ such that the
ring extension $R\subset \cO(\Phi^{-1}(U))$ is generically separable.
For a finite group $G$, $\Phi$ is said to be \emph{$G$-cover} if in
addition there exists a group homomorphism $G\rightarrow \Aut_X(Y)$
which acts transitively on the geometric generic fibers of $\Phi$.

For an integral scheme $X$, $\pi_1(X)$ will denote the algebraic
fundamental group of $X$ with respect to the generic point. For an
affine variety $X$ over a field $k$, $k[X]$ will denote the coordinate
ring of $X$ and $k(X)$ will denote the function field. For a scheme
$X$ and a point $x \in X$, let $\hat{\cK}_{X,x}$ denote the fraction
field of complete local ring $\hat{\cO}_{X,x}$ whenever the latter is
a domain.  For domains $A\subset B$, $\overline{A}^B$ will denote the
integral closure of $A$ in $B$.

For a scheme $X$, let $\cP(X)$ denote the category of coherent
sheaves of $\cO_X$-modules, $\cAP(X)$ denote the category of
coherent sheaves of $\cO_X$-algebras and $\cSP(X)$ denote
the subcategory of $\cAP(X)$ for which the sheaves are generically
separable and locally free.  For a finite group $G$, let $G\cP(X)$ denote the category
of generically separable coherent locally free sheaves of
$\cO_X$-algebras $S$ together with a $G$-action which is transitive on
the geometric generic fibers of $\spec_{\cO_X}(S) \rightarrow X$. For
a ring $R$, we may use $\cP(R)$ instead of $\cP(\spec(R))$, etc. Given
categories $\cA$, $\cB$ and $\cC$ and functors from
$\cF:\cA\rightarrow\cC$ and $\cG:\cB \rightarrow \cC$, we define the
fiber product category $\cA\times_{\cC}\cB$ whose objects are triples
$(A,B,C)$, where $A$, $B$ and $C$ are objects of $\cA$, $\cB$ and
$\cC$ respectively together with isomorphisms of $C$ with $\cF(A)$ and
with $\cG(B)$ in $\cC$, morphisms are triples $(a,b,c)$, where $a$,
$b$ and $c$ are morphisms in $\cA$, $\cB$ and $\cC$ respectively, so
that $\cF(a)$ and $\cG(b)$ under the functors $\cF$ and $\cG$ are
morphism in $\cC$ which agrees with $c$ in the natural way. That is,
the following two squares commute.
$$
\xymatrix{
  C\ar[d] \ar[r]^c       & C'\ar[d]   & & & C\ar[d] \ar[r]^c     & C'\ar[d]\\
  \cF(A) \ar[r]_{\cF(a)} & \cF(A') & & & \cG(B) \ar[r]_{\cG(b)} &
  \cG(B') }
$$

\section{abelianization}

Let $X=\spec(A)$ be a normal affine algebraic variety over an
algebraically closed field $k$ of characteristic $p > 0$. For a ring
$R$, by $(W_n(R), +, .)$ we denote the ring of Witt vectors of length
$n$ over $R$. Let $F$ be the Frobenius map on $W_n(R)$ which sends
$(a_1,..,a_n)$ to $(a_1^p,..,a_n^p)$ and $P:W_n(R) \rightarrow W_n(R)$
be the map of abelian groups sending $(a_1,..,a_n) \mapsto
(a_1^p,..,a_n^p) - (a_1,..a_n)$, i.e., $F$ - Identity. For a detailed
account of Witt vectors, readers are advised to see [Jac, Chapter 8].
In fact, some of the ideas for the proof of Lemma 3.3 below come from
this source.

Let $G=\pi_1 ^{ab}(X)$ and let $G_p$ be the maximal $p$ quotient of
$G$, i.e., $G_p$ is the quotient group of $G$ such that every finite
$p$-group quotient of $G$ factors through $G_p$. Note that $G_p$ is
also a subgroup of $G$.  Let $G_W=\Hom(\dirlim W_n(A)/ P(W_n(A)),
S^1)$, where the group homomorphism from $W_n(A)/P(W_n(A))$ to
$W_{n+1}(A)/P(W_{n+1}(A))$ is given by $[(a_1,..,a_n)] \mapsto
[(0,a_1,..a_n)]$.

\begin{thm}
  The p-part of the abelianization of the fundamental group $G_p$ is
  isomorphic to $G_W$.
\end{thm}

First we need a few lemmas.  For a sheaf of rings $\cF$ of
characteristic $p$ on a topological space X, let $W_n(\cF)$ denote the
sheaf which assigns to an open set $U$, the ring $W_n(\cF (U))$. A
version of the following lemma can be found in [Se1].

\begin{lm}(\text{\bf Serre}) Let $A$ be a noetherian ring of characteristic $p$. 
  Let $B$, also of characteristic $p$, be a finite ring extension of
  $A$ then for every $n\ge 1$, $H_{et}^k(\spec(A),W_n(\cB))=0$, $k >
  0$, where $\cB=\theta_* \cO_{\spec(B)}$ and $\theta:\spec(B)
  \rightarrow spec(A)$ is the morphism induced from $A \into B$.
\end{lm}

\proof $B$ is a finite module over $A$, hence $\cB$ is coherent sheaf
over $\spec(A)$. We shall use induction on $n$ to prove the lemma.
Note that $W_1(\cB)=\cB$.  By the Serre's vanishing theorem and the
fact that \'etale cohomology of coherent sheaves agrees with the
Zariski cohomology (see [Mil, 3.7, 3.8, page 114]), the lemma holds
for $n=1$. For the induction step, consider the following exact
sequence.
$$0 \longrightarrow \cB \longrightarrow W_{n+1}(\cB) \longrightarrow
W_n(\cB) \longrightarrow 0$$
where for a fixed open set $U$ the
surjection (on the level of rings) is given by $(b_1,..,b_n,b_{n+1})
\mapsto (b_1,..,b_n)$, clearly the kernel is $\{(0,..,0,b) \in
W_{n+1}(\cB(U)): b\in \cB(U)\} \cong \cB(U)$.  This induces a long
exact sequence
$$... \rightarrow H_{et}^k(\spec(A),\cB) \longrightarrow
H_{et}^k(\spec(A),W_{n+1}(\cB)) \rightarrow
H_{et}^k(\spec(A),W_n(\cB)) \longrightarrow ...$$
By the induction
hypothesis on $n$, $H_{et}^k(\spec(A),W_n(\cB))=0$ for $k>0$ and hence
$H_{et}^k(\spec(A),W_{n+1}(\cB))=0$, for all $k>0$.  \qed

Now extending the Artin-Schrier theory to the Witt vectors we get the
following result.

\begin{lm} Let $A$ be a finitely generated normal domain over an 
  algebraically closed field $k$ of characteristic $p > 0$.  Let
  $\pi_1(X)$ be the fundamental group of $X=\spec(A)$ and
  $(W_n(A),+,.)$ be the ring of Witt vectors of length $n$.  Let $P$
  be the additive group endomorphism, $F - Id$, of $W_n(A)$. Then for
  every $n \ge 1$, we have a natural isomorphism $W_n(A)/P(W_n(A))
  \stackrel{\Phi}{\longrightarrow} \Hom(\pi_1(X),W_n(\ZZ/p\ZZ))$ so
  that the following diagram commutes.
  
  $$
  \xymatrix{
    W_n(A)/P(W_n(A)) \ar[r] \ar[d] \ar@{}[dr]|\chash & \Hom(\pi_1(X),W_n(\ZZ/p\ZZ)) \ar[d]\\
    W_{n+1}(A)/P(W_{n+1}(A)) \ar[r] & \Hom(\pi_1(X),W_{n+1}(\ZZ/p\ZZ))
  }$$
  where the first vertical map sends $[(a_1,..,a_n)]$ to
  $[(0,a_1,..,a_n)]$ and second vertical map is induced by inclusion
  of $W_n(\ZZ/p\ZZ)$ in $W_{n+1}(\ZZ/p\ZZ)$
\end{lm}
\proof Let $K$ be an algebraic closure of the fraction field $\QF(A)$
of $A$. Let $K^{un}$ be the compositum of all subfields of $K$, $L$,
so that $L/\QF(A)$ is finite field extension and the integral closure
in $L$, $\overline{A}^L$ is an unramified (hence \'etale) ring
extension of $A$.  We define a map $\phi$ from $W_n(A) \rightarrow
\Hom(\pi_1(X),W_n(\ZZ/p\ZZ))$. Given $(a_1,..,a_n) \in W_n(A)$, let
$(r_1,..,r_n) \in W_n(K)$ be such that $P(r_1,..,r_n)=(a_1,..,a_n)$.
Note that $r_1,..,r_n \in K^{un}$; to see this first observe that $r_i
^p - r_i \in A[r_1,..,r_{i-1}]$ for each $i\ge 1$. We know that the
ring extension given by the polynomial of the form $Z^p - Z - a$, for
$a\in A$ is unramified over $A$.  $\pi_1(X)=\Gal(K^{un}/\QF(A))$, so
$\pi_1(X)$ acts on $K^{un}$ fixing $\QF(A)$.  For $g\in \pi_1(X)$, let
$\phi(a_1,..,a_n)(g)=(gr_1,..,gr_n)-(r_1,..,r_n)$. Note that
\begin{align*}
  P((gr_1,..,gr_n)-(r_1,..,r_n)) &= P((gr_1,..,gr_n))-P(r_1,..,r_n)\\
  &=gP(r_1,..,r_n)-(a_1,..,a_n)=0. \\
\end{align*}
Since $(r_1 ^p,..,r_n ^p)-(r_1,..,r_n)$ is given by polynomial in
$r_1,..,r_n$ with integer coefficients.  Hence
$F((gr_1,..,gr_n)-(r_1,..,r_n))=(gr_1,..,gr_n)-(r_1,..,r_n)$, which
gives us $(gr_1,..,gr_n)-(r_1,..,r_n) \in W_n(\ZZ/p\ZZ)$. To see that
$\phi(a_1,..,a_n)$ is independent of the choice of $(r_1,..,r_n)$, let
$(s_1,..,s_n)$ be such that $P(s_1,..,s_n)=(a_1,..,a_n)$ then
$(r_1,..,r_n)-(s_1,..,s_n) \in W_n(\ZZ/p\ZZ)$, hence fixed by $g$. So
$g((r_1,..,r_n)-(s_1,..,s_n))=(r_1,..,r_n)- (s_1,..,s_n)$ which yields
$(gr_1,..,gr_n)-(r_1,..,r_n)=(gs_1,..gs_n)-(s_1,..s_n)$. Next we shall
see $\phi(a_1,..,a_n)$ is a homomorphism from $\pi_1(X)$ to
$W_n(\ZZ/p\ZZ)$.  Let $g,h \in \pi_1(X)$ then
\begin{align*}
  \phi(a_1,..,a_n)(gh)&=(ghr_1,..,ghr_n)-(r_1,..,r_n) \\
  &=(ghr_1,..,ghr_n)-(hr_1,..,hr_n)+ (hr_1,..,hr_n)-(r_1,..r_n) \\
  &=\phi(a_1,..,a_n)(g)+\phi(a_1,..,a_n)(h)
\end{align*}
since $P(hr_1,..,hr_n)=(a_1,..,a_n)$. Now we shall see $\phi$ is a
homomorphism. To simplify notation, we may write $\underline{a}$ for
$(a_1,..,a_n)$.  Let $\underline{a}, \underline{b} \in W_n(A)$ and
$\underline{r}, \underline{s} \in W_n(K)$ be such that
$P(\underline{r})=\underline{a}, P(\underline{s})=\underline{b}$ then
$P(\underline{r}+\underline{s})=\underline{a}+\underline{b}$. Hence
$\phi(\underline{a}+\underline{b})=\phi(\underline{a})+\phi(\underline{b})$.
To determine the kernel of $\phi$, note that $\phi(\underline{a})=0$
iff $g\underline{r}=\underline{r}$ for all $g \in G$, i.e.
$\underline{r} \in W_n(\QF(A))$ but $A$ is normal, hence
$\underline{r} \in W_n(A)$. Hence the kernel of $\phi$ is $P(W_n(A))$.
Let $\Phi$ be the induced map on $W_n(A)/P(W_n(A))$. The fact that the
diagram commutes is obvious by the construction.  So to complete the
proof it suffices to show that $\phi$ is surjective.

Let $\alpha \in \Hom(\pi_1(X),W_n(\ZZ/p\ZZ))$, we shall find a Witt
vector $(a_1,..,a_n)$ so that $\alpha=\phi(a_1,..,a_n)$. Note that
$\alpha$ corresponds to a Galois \'etale extension $B$ of $A$ with
Galois group of $\QF(B)$ over $\QF(A)$ being $\image(\alpha)(=H \text{
  say})$. Let $SW_n(B)=\{(r_1,..,r_n) \in W_n(B):P(r_1,..,r_n) \in
W_n(A)\}$.  Clearly $W_n(A) \into SW_n(B)$.  Let $\hat{H}=Hom(H,S^1)$
be the character group of $H$. For $\underline{r} \in SW_n(B)$ and $h
\in H$ define $\chi_{\underline{r}} (h)=h\underline{r} - \underline{r}
:= g\underline{r} - \underline{r}$ where $g$ is any element of $\alpha
^{-1}(h)$. As noted earlier $\chi_{\underline{r}}$ is a character
(after identifying $W_n(\ZZ/p\ZZ)$ with the unique cyclic subgroup of
$S^1$). We know that $\hat{H} \cong H$. Also as seen earlier $\lambda:
\underline{r} \mapsto \chi_{\underline{r}}$ is a group homomorphism
from $SW_n(B)$ to $\hat{H}$ whose kernel is precisely $W_n(A)$. So we
have an exact sequence
$$0 \longrightarrow W_n(A) \longrightarrow SW_n(B) \longrightarrow
\hat{H}$$
where the last homomorphism is $\lambda$. Next we shall show
that $\lambda$ is surjective. Since $\hat{H} \cong H$ is a quotient of
$\pi_1(X)$, we have $H^1(\hat{H},W_n(B)) \into H^1(\pi_1(X),W_n(B))$
[Wei, 6.8.3]. By the Hochschild-Serre spectral sequence
$H^1(\pi_1(X),W_n(B))=H^1(\pi_1(X),H_{et}^0(X,W_n(\cB)))$ which embeds
into $H_{et}^1(X,W_n(\cB))$ (see [Mil, 2.21(b), page 106]) and by
previous lemma, $H_{et}^1(X,W_n(\cB))=0$. So $H^1(\hat{H},W_n(B))=0$,
i.e., every cocycle is a coboundary. And if $\chi \in \hat{H}$ then
$\chi(hh')=\chi(h)+\chi(h')=h\chi(h')+\chi(h)$, i.e., it is a cocycle
hence a coboundary.  So there exists $\underline{r} \in W_n(B)$ such
that $\chi(h)=h\underline{r}-\underline{r}$.  Since $\chi(h) \in
W_n(\ZZ/p\ZZ)$, $P(\chi(h))=0$, $\forall h\in H$, i.e.
$hP(\underline{r})=P(\underline{r}), \forall h \in H$. Since $A$ is
normal this means $P(\underline{r}) \in W_n(A)$, hence $\underline{r}
\in SW_n(B)$. This proves $\lambda$ is surjective.  So we have
$SW_n(B)/W_n(A) \cong \hat{H} \cong H$. Since $H$ is a subgroup of
$W_n(\ZZ/p\ZZ)$, $H$ has a generator of the type $h=(0,..0,1,0,..,0)$.
Let the coset $(r_1,..,r_n) +W_n(A)$ be a generator of
$SW_n(B)/W_n(A)$. It follows that $\chi_{\underline{r}}$ is a
generator of $\hat{H}$ and hence $h_1:=\chi_{\underline{r}}(h)$ is a
generator of $H$. So there is an $h_2 \in W_n(\ZZ/p\ZZ)$ such that
$h_1.h_2=h$. Let $g \in\alpha^{-1}(h)$ then
$g\underline{r}-\underline{r} = \chi_{\underline{r}}(h)=h_1$.  Let
$(a_1,..,a_n)=P(h_2.\underline{r})$. This Witt vector is our candidate
for preimage of $\alpha$, we shall show $\alpha=\phi(a_1,..,a_n)$.  By
assumption $\alpha(g)=h$ and
\begin{align*}
  \phi(a_1,..,a_n)(g)&=g(h_2.\underline{r})-h_2.\underline{r}\\
  &=gh_2.g\underline{r}-h_2.\underline{r} \\
  &=h_2.g\underline{r}-h_2.\underline{r} \\
  &=h_2.(g\underline{r}-\underline{r}) \\
  &=h_2.h_1=h
\end{align*}
For arbitrary $g_1 \in G$, $\alpha(g_1)=h+..+h$ say $k$ times, since
$H$ is cyclic. Then
\begin{align*}
  \phi(a_1,..,a_n)(g_1)&=h_2.(g_1\underline{r}-\underline{r}) \\
  &=h_2.(\alpha(g_1)\underline{r}-\underline{r}) \\
  &=h_2.(\chi_{\underline{r}}(h+..+h)) \\
  &=h_2.(\chi_{\underline{r}}(h)+..+\chi_{\underline{r}}(h)) \\
  &=h_2.(h_1+..+h_1)=h+..+h
\end{align*}
So $\phi(a_1,..,a_n)$ agrees with $\alpha$ on whole of $\pi_1(X)$.
\qed

\proof{\bf (Theorem 3.1)} We know that $G_p \cong \Hom( \Hom_{cont}
(G_p,S^1),S^1)$ by Pontriagan duality.  Let $K_n$ be the compositum of
function fields of finite Galois \'etale extension of $\spec(A)$ with
Galois group a subgroup of $(\ZZ/p^n\ZZ)^m$ for some m. And let $G_n$
be the Galois group of $K_n$ over $\QF(A)$.  The natural group
homomorphism from $G_{n+1}$ to $G_n$ corresponding to the Galois
extension $K_{n+1} \supset K_n \supset \QF(A)$ makes $(G_n)_{n \ge 1}$
into an inverse system and $G_p=\invlim G_n$. So we have
$\Hom_{cont}(G_p,S^1) \cong \Hom_{cont}(\invlim G_n,S^1)$. Since
$\Hom$ is contravariant and dual of inverse limit is direct limit,
this is isomorphic to $\dirlim$$\Hom_{cont}(G_n,S^1)$ Now since $G_n$
is a $p^n$ torsion group and $W_n(\ZZ/p\ZZ)$ and be identified as the
unique cyclic subgroup of $S^1$ of order $p^n$, we get
$\Hom_{cont}(G_n,S^1) \cong \Hom_{cont}(G_n,W_n(\ZZ/p\ZZ))$. Also
$\Hom_{cont}(G_n,W_n(\ZZ/p\ZZ))\cong$
$\Hom_{cont}(G_p,W_n(\ZZ/p\ZZ))$, since all maps from $G_p$ to
$W_n(\ZZ/p\ZZ)$ factors through $G_n$.  Similarly, $\Hom_{cont}(G_p,
W_n(\ZZ/p\ZZ))$ is isomorphic to $\Hom_{cont}(\pi_1(X),
W_n(\ZZ/p\ZZ))$, since $W_n(\ZZ/p\ZZ)$ is an abelian $p$-group, all
homomorphisms from $\pi_1(X)$ to $W_n(\ZZ/p\ZZ)$ factors through
$G_p$. Now by the previous lemma $\dirlim$ $\Hom_{cont}(\pi_1(X),
W_n(\ZZ/p\ZZ))$ is isomorphic to $\dirlim W_n(A)/P(W_n(A))$. \qed

The following result is a corollary of a classical result of
Grothendieck [SGAI, XIII, Corollary 2.12, page 392].

\begin{thm}(\text{\bf Grothendieck})
  The prime to $p$ part of the abelianization of the fundamental group
  of an affine curve $C=\spec(A)$ over an algebraically closed field
  $k$ of characteristic $p>0$ is given by $\displaystyle
  \bigoplus_{i=1}^{2g+r-1} (\bigoplus_{l\ne p \text{ prime}}\ZZ_l)$
  where $g$ is the genus of the smooth compactification curve and $r$
  is the number of points in the compatification which are not in $C$.
\end{thm}

\begin{cor}
  Under the assumption of the previous theorem, abelianization of the
  fundamental group of $C$, $\pi_1^{ab}(C)$, is given by $$
  \Hom(\dirlim W_n(A)/ P(W_n(A)), S^1) \bigoplus \displaystyle
  \bigoplus_{i=1}^{2g+r-1} (\bigoplus_{l\ne p \text{ prime}}\ZZ_l)$$
\end{cor}
\proof This follows directly from Theorem 3.1 and Theorem 3.4.

\begin{cor}
  Since the rank of $\pi_1^{ab}(C)$ is same as the cardinality of $k$,
  we get another proof of a known result that $\pi_1^{ab}(C)$
  determines the cardinality of the base field. In fact just the
  $p$-part determines $W_n(A)/P(W_n(A))$ for all $n$.
\end{cor}
\proof This is a direct consequence of Lemma 2.3.

\section{Group theory}

Let $k$ be an algebraically closed field of characteristic
$p > 0$. Let $\pi_1(C)$ be the algebraic fundamental group of a smooth
affine curve $C$ over $k$ and $\pi_1 ^c(C)=[\pi_1(C),\pi_1(C)]$ be the
commutator subgroup.

\begin{thm}(Main theorem) Let $C$ be an irreducible smooth affine curve over a countable 
  algebraically closed field of characteristic $p$. Then $\pi_1 ^c(C)$ is 
  free of countable rank.
\end{thm}

This will be reduced to finding solutions of certain kinds of
embedding problems.  Before that we need a group theoretic result
which connects ``freeness'' of a profinite group to solving embedding
problems. Below are certain results in this direction.

\begin{thm} (\text{\bf Iwasawa [Iwa, p.567], [FJ, Corollary 24.2]})
  A profinite group $\pi$ of countably infinite rank is free if and
  only if every finite embedding problem for $\pi$ has a proper
  solution.
\end{thm}

This was generalised by Melnikov and Chatzidakis for any cardinality (cf [Jar, Theorem 2.1]).
The Melnikov-Chatzidakis result says that for an infinite cardinal $m$, a profinite group 
$\pi$ is free of rank $m$ if and only if every finite nontrivial embedding problem for 
$\pi$ has exactly $m$ solution. Following is a variant of this result which has been 
proved in [HS].

\begin{thm} (\text{\bf [HS, Theorem 2.1]})
Let $\pi$ be a profinite group and let $m$ be an infinite cardinal. Then $\pi$ is a
free profinite group of rank $m$ if and only if the following conditions are satisfied:\\
(i) $\pi$ is projective.\\
(ii) Every split embedding problem for $\pi$ has exactly $m$ solution.\\
\end{thm}

We shall anyway see a standard argument which reduces the problem of finding proper
solutions of an embedding problem for a projective profinite group to finding proper
solutions of a \emph{split} embedding problem with the same kernel.

Let $C$ be a smooth affine curve over an algebraically closed field $k$ of cardinality 
$m$. Since $k(C)$, the function field of $C$, is also of cardinality $m$, there are
only $m$ polynomials over $k(C)$.  Hence the absolute
Galois group of $k(C)$ is the inverse limit of finite groups over a
set of cardinality $m$ and hence has generating set of cardinality $m$
(generating set in the topological sense).
So $\pi_1(C)$, being a quotient of the absolute Galois group of
$k(C)$, is $m$ genrated and hence $\pi_1 ^c(C)$ is $m$
generated. So to prove that $\pi_1 ^c(C)$ is free, it suffices to show
that every embedding problem for $\pi_1 ^c(C)$ has $m$ proper solution (and 
just one solution if $m$ is the countable cardinal),
since this implies $\pi_1 ^c(C)$ has rank exactly $m$.
So given an embedding problem:
$$(*) \xymatrix{
  &          &                        &\pi_1 ^c(C) \ar @{-->}[dl]_{\psi} \ar[d]^{\phi}\\
  1\ar[r] & H \ar[r] & \Gamma \ar[r]_{\alpha} & G \ar[r] \ar[d] & 1\\
  & & & 1 }
$$
we need to find $\card(k)=m$ proper solution (and just one solution if $m$ 
is the countable cardinal) for every finite group $G$, $\Gamma$ and $H$.

Before that, we shall show that $\pi_1 ^c(C)$ is projective and use it
to reduce to the case where $(*)$ splits.

\begin{pro}
  For an irreducible smooth affine curve $C$ over $k$,
  $\pi_1 ^c(C)=[\pi_1(C), \pi_1(C)]$ is a projective group. More
  explicitly, given
  $$\xymatrix{
    &\pi_1 ^c(C) \ar @{-->}[dl]_{\exists \psi} \ar @{->>}[d]^{\phi} \\
    \Gamma \ar @{->>}[r]_{\alpha} & G }$$
  surjections $\phi$ and
  $\alpha$ to a finite group $G$ from $\pi_1 ^c(C)$ and another finite
  group $\Gamma$ respectively, there exist a group homomorphism $\psi$
  from $\pi_1 ^c(C)$ to $\Gamma$ so that the above diagram commutes,
  i.e., $\alpha \circ \psi =\phi$
\end{pro}
\proof Let $K^{ab}$ be the compositum of the function fields of abelian
\'etale covers of $C$, i.e., the compositum of all $L$, $k(C)\subset
L\subset K^{ab}$ with $L/k(C)$ finite, the integral closure of $k[C]$
in $L$, $\overline{k[C]}^L$, is \'etale extension of $k[C]$ and
$\Gal(L/k(C))$ is abelian.

A surjection $\phi:\pi_1 ^c(C)\rightarrow G$ corresponds to a Galois
field extension $M/K^{ab}$ with the Galois group $\Gal(M/K^{ab})=G$
and $M\subset K^{un}$, where $K^{un}$ is the compositum of the function
fields of all \'etale covers of $C$.

Since $M/K^{ab}$ is a finite field extension, there exist $L$,
$k(C)\subset L\subset K^{ab}$, a finite Galois extension of $k(C)$ and
$L'$, a Galois extension of $L$ with $\Gal(L'/L)=G$ and $L'K^{ab}=M$.
Let $\pi_1 ^L=\pi_1(\spec(\overline{k[C]}^L))$.

So we have the following tower of fields.\\
\centerline{
  \xymatrix{K^{un}                       & M=L'K^{ab} \ar[l] \\
    K^{ab}\ar[u]_{\pi_1 ^c(C)} \ar[ru]_G & L' \ar[u]\\
    L \ar@/^2pc/[uu]^{\pi_1 ^L} \ar[u] \ar[ru]_G\\
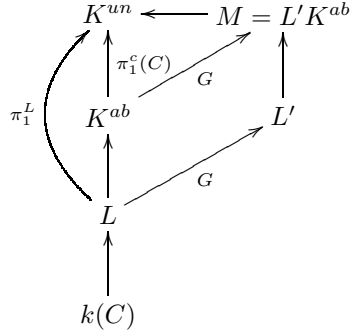
    k(C)\ar[u] } } \centerline{Fig. 1} Moreover
$\Gal(K^{un}/K^{ab})=\pi_1 ^c(C)$, $\Gal(K^{un}/L)=\pi_1 ^L$ and
$\pi_1 ^c(C)$ is a subgroup of $\pi_1 ^L$.  The field extension $L'/L$
gives a surjection $\tilde{\phi}:\pi_1 ^L \rightarrow G$. Since $L'/L$
is a descent of the field extension $M/K^{ab}$, $\tilde{\phi}|_{\pi_1
  ^c(C)}=\phi$.  By [Se2, Proposition 1], which says the fundamental
group of any affine curve is projective, we have $\pi_1
^L:=\pi_1(\spec(\bar{k[C]}^L))$ is projective. So there exists a lift,
$\tilde{\psi}$, to
$\Gamma$ of $\tilde{\phi}$. i.e.,\\
\centerline{ \xymatrix{
    &\pi_1 ^L \ar @{-->}[dl]_{\exists \tilde{\psi}} \ar @{->>}[d]^{\tilde{\phi}}\\
    \Gamma \ar @{->>}[r]_{\alpha} & G } } with
$\alpha\circ\tilde{\psi}=\tilde{\phi}$. So
$\alpha\circ\tilde{\psi}|_{\pi_1 ^c(C)}=\tilde{\phi}|_{\pi_1
  ^c(C)}=\phi$. So $\tilde{\psi}|_{\pi_1 ^c(C)}$ gives a lift of
$\phi$.  \qed

To reduce to the case where $(*)$ splits, let $G'=\im{\psi}$, where
$\psi$ is as in Propostion 4.3.  $G'$ acts on $H$ by conjugation,
since $H$ is a normal subgroup of $\Gamma$. Let $\Gamma'=H\rtimes G'$
then we have a natural surjection $\beta:\Gamma'\rightarrow\Gamma$
given by $(h,g) \mapsto hg$. So if we have a proper solution $\theta'$
for the embedding
problem, \\
\centerline{ \xymatrix{
    &          &               &\pi_1 ^c(C) \ar @{-->}[dl]_{\theta '} \ar[d]^{\tilde{\psi}}\\
    1\ar[r] & H \ar[r] & \Gamma' \ar[r] & G' \ar[r]\ar[d] & 1\\
    & & & 1 } } then $\theta=\beta\circ\theta'$ provides a proper
solution to $(*)$. Note that this reduction holds for any projective profinite group
not necessarily $\pi_1 ^c(C)$.

From now onwards we shall assume that all our embedding problems are
split embedding problems unless otherwise stated. Our proof inducts on the cardinality
of $H$. So whenever we encounter an embedding problem which may not be split, without 
explicitly stating we shall assume that the embedding problem has been replaced by an
appropriate split embedding problem with the same kernel.

\begin{thm}
Let $\pi$ be any projective profinite group of rank exactly $m$, then $\pi$ is free of 
rank $m$ if and only if for any finite group $\Gamma$ and any minimal normal sugroup $H$ 
of $\Gamma$, the embedding problem
$$\xymatrix{
  &          &                        &\pi \ar @{-->}[dl]_{\psi} \ar[d]^{\phi}\\
  1\ar[r] & H \ar[r] & \Gamma \ar[r]_{\alpha} & G \ar[r] \ar[d] & 1\\
  & & & 1 }
$$
has $m$ distinct solutions (and atleast one solution if $m$ 
is the countable cardinal) in the following three cases:\\
$(1)$ $H$ is a quasi-$p$ perfect group, i.e. $H=[H,H]$.\\
$(2)$ $H$ is an abelian $p$-group.\\
$(3)$ $H$ is a prime-to-$p$ group.\\
\end{thm}

\proof In view of Theorem 4.3 (and Theorem 4.2 if $m$ is the countable cardinal),
``only if part'' is trivial and for ``if part''  it is enough to 
show the embedding problem for $\pi$ has $m$ distinct proper solutions any finite group $H$.
We induct on the cardinality of $H$. Suppose $H$ is not minimal normal subgroup. 
Let $H_1$ be a proper nontrivial subgroup of $H$ and $H_1$ is a normal subgroup of $\Gamma$.
Then we have the following two proper nontrivial embeddding problems.
$$
\xymatrix{
  &               &                    &\pi \ar @{-->}[dl] \ar[d] \\
  1\ar[r] & H/H_1 \ar[r] & \Gamma/H_1 \ar[r] & G \ar[r]\ar[d] & 1 \\
  & & & 1 }$$
and
$$
\xymatrix{
  &             &               &\pi \ar @{-->}[dl] \ar[d] \\
  1\ar[r] & H_1 \ar[r] & \Gamma \ar[r] & \Gamma/H_1 \ar[r]\ar[d] & 1 \\
  & & & 1 }$$

Since the cardinality of $H_1$ and $H/H_1$ is strictly smaller than the 
cardinality of $H$, by induction hypothesis, we have $m$ distinct proper 
solutions to these embedding problems. Hence we have $m$ distinct proper
solutions to the embedding problem.
Hence we may assume $H$ is a nontrivial  minimal normal subgroup of $\Gamma$. 
So $H\cong \SSS\times..\times\SSS$ for some finite simple group $\SSS$.
If $\SSS$ is prime-to-$p$ then $H$ is prime-to-$p$, hence we are done by case $(3)$.
If $\SSS$ is quasi-$p$ nonabelian group then $H$ being product of perfect groups is perfect. 
So we are done by case $(1)$. And finally if $\SSS$ is quasi-$p$ abelian then $\SSS\cong \ZZ/p\ZZ$. Hence $H$ is abelian $p$-group and we are done by case $(2)$.
\qed

\proof{\bf (Theorem 4.1)}
In view of Theorem 4.4, Theorem 4.1 follows from the previous theorem, 
once we show that there exist a solution to the embedding problem for $\pi_1 ^c(C)$ for
the case $(1)$, $(2)$ and $(3)$
of the previous theorem. The case $(1)$ and $(2)$ will be proved in Section 5 and
$(3)$ will be proved in Section 6.
\qed

Before proving that a solution to embedding problem $(*)$ exists in
the above three cases, we shall prove the following group theory
result. This will be used in the next section.
\begin{lm}
  Given any finite abelian $p$-group $A$ there exists a finite
  $p$-group $B$ such that the commutator of $B$, $[B,B]=A$.
\end{lm}

\proof Since $A$ is an abelian $p$-group. $A$ is a direct sum of
cyclic $p$-groups. Observe that the commutator of the group $B_1\times
B_2$ is isomorphic to $[B_1, B_1]\times [B_2, B_2]$ for any two groups
$B_1$ and $B_2$.  So we may assume $A$ is a cyclic $p$-group, say
$\ZZ/p^m\ZZ$. Consider the Heisenberg group over $\ZZ/p^m\ZZ$, i.e.,
the group of $3\times 3$ upper triangular matrices with diagonal
entries $1$. It is a group of order $p^{3m}$ generated by the matrices
\[
\begin{pmatrix}
  1 & 1 & 0\\
  0 & 1 & 0\\
  0 & 0 & 1
\end{pmatrix},\ \ \begin{pmatrix}
  1 & 0 & 0\\
  0 & 1 & 1\\
  0 & 0 & 1
\end{pmatrix} \]
and one could easily check that the commutator of this group is the
subgroup generated by \[
\begin{pmatrix}
  1 & 0 & 1 \\
  0 & 1 & 0 \\
  0 & 0 & 1
\end{pmatrix} \]
which is clearly isomorphic to $(\ZZ/p^m\ZZ,+)$.  \qed

The construction of such a group using Heisenberg matrices was pointed
out to me by a friend Sandeep Varma and also by Prof. Donu Arapura.

\section{Quasi-p embedding problem}

In this section we show that the split embedding problem has a
solution if $H$ is a perfect quasi-$p$ group or $H$ is a $p$-group. We
shall begin by stating some results on quasi-$p$ embedding problems.

\begin{thm}({\bf Florian Pop, [Pop], [Ha3, Theorem 5.3.4], [Ha6, Corollary 4.6]})
  Let $k$ be an algebraically closed field of characteristic $p > 0$, $\card(k)=m$,
  and let $C$ be an an irreducible affine smooth curve over $k$. Then
  every quasi-$p$ embedding problem for $\pi_1(C)$ has $m$ distinct proper
  solutions.
\end{thm}

\begin{thm}({\bf [Ha5, Theorem 1b]})
  Let $\pi$ be a profinite group such that $H^1(\pi,P)$ is infinite
  for every finite elementary abelian $p$-group $P$ with continuous
  $\pi$-action. Then every $p$-embedding problem for $\pi$ has a
  proper solution if and only if every $p$-embedding problem has a
  weak solution (equivalently, $p$ cohomological dimension of $\pi$,
  $\cd_p(\pi)\le 1$).
\end{thm}

\begin{thm}
  The following split embedding problem has $\card(k)=m$ proper solutions
  $$
  \xymatrix{
    &          &               &\pi_1 ^c(C) \ar[d] \ar @{-->} [dl]\\
    1\ar[r] & H \ar[r] & \Gamma \ar[r] & G \ar[r]\ar[d] & 1\\
    & & & 1 }
  $$
  Here $H$ is a quasi-$p$ perfect group (i.e. $[H,H]=H$) and $\pi_1
  ^c(C)$ is the commutator of the algebraic fundamental group of an
  irreducible smooth affine curve $C$ over an algebraically closed
  field $k$ of characteristic $p$.
\end{thm}

\proof As in Proposition 4.4 (also see Fig. 1), let $K^{un}$ denote
the compositum (in some fixed algebraic closure of $k(C)$) of the function
fields of all \'etale Galois covers of $C$.  And let $K^{ab}$ be the
subfield of $K^{un}$ obtained by considering only abelian \'etale
covers of $C$. In terms of Gaolis theory, $\pi_1 ^c(C)$ is
$\Gal(K^{un}/K^{ab})$. So giving a surjection from $\pi_1 ^c(C)$ to
$G$ is same as giving a Galois extension $M\subset K^{un}$ of $K^{ab}$
with Galois group $G$. Since $K^{ab}$ is an algebraic extension of
$k(C)$ and $M$ is a finite extension of $K^{ab}$, we can find a finite
abelian extension $L \subset K^{ab}$ of $k(C)$ and $L'\subset K^{un}$
a Galois extension of $L$ with Galois group $G$ so that $M=K^{ab}L'$. Let $X$ be the
normalization of $C$ in $L$ and $ \Phi_X$ be the normalization
morphism. Then $X$ is an \'etale abelian cover of $C$ and the function
field, $k(X)$, of $X$ is $L$. Let $W_X$ be the normalization of $X$ in
$L'$ and $\Psi_X$ corresponding normalization morphism. Then $\Psi_X$
is \'etale and $k(W_X)=L'$.

By Theorem 5.1 applied to the affine curve $X$ and translating the
conclusion into Galois theory, we conclude that there exist $m$ distinct
smooth irreducible \'etale $\Gamma$-covers. Each one of these $\Gamma$-cover, 
$Z$, of $X$ is such that $Z/H=W_X$. Clearly
$k(Z)\subset K^{un}$.  
We also have $\Gal(k(Z)K^{ab}/K^{ab})\subset\Gamma$ and by assumption
$\Gal(k(W_X)K^{ab}/K^{ab})=G$.
Moreiver, Galois
group of $k(Z)/k(W_X)$ is $H$ which is a perfect group and
$k(W_X)K^{ab}/k(W_X)$ is a pro-abelian extension. Hence they are
linearly disjoint, so $\Gal(k(Z)K^{ab}/k(W_X)K^{ab})=H$. 
So $\Gal(k(Z)K^{ab}/K^{ab})=\Gamma$. Also if $Z$ and $Z'$ are two distinct
solutions then the $\Gal(k(Z)k(Z')/k(W_X))$ is quotient of $H\times H$ 
and hence perfect. So $\Gal(k(Z)k(Z')K^{ab}/k(W_X)K^{ab} = 
\Gal(k(Z)k(Z')/k(W_X))$ and consequently $k(Z)K^{ab}$ and $k(Z')K^{ab}$
are distinct fields. 
\qed

Now we shall consider the case when $H$ is an abelian $p$-group. 

\begin{lm}
  Let $P$ be any nonzero finite abelian $p$-group, then there exist $m$ 
  distict surjections from $\pi_1 ^c(C)$ to $P$.
\end{lm}

\proof Let $n \ge 1$ be a natural number. By Lemma 4.6 there exist a
$p$-group $P_1$ such that its commutator $[P_1,P_1]=P$. By Theorem
5.1 there exist $m$ distinct surjective homomorphisms from $\pi_1(C)$
to $P_1$ and clearly the commutator $\pi_1 ^c(C)$ surjects onto
$[P_1,P_1]=P$ under their restrictions. 

\begin{thm}
  The following split embedding problem has $\card(k)=m$ proper solutions
  $$
  \xymatrix{
    &          &               &\pi_1 ^c(C) \ar[d] \ar @{-->} [dl]\\
    1\ar[r] & H \ar[r] & \Gamma \ar[r] & G \ar[r]\ar[d] & 1\\
    & & & 1 }
  $$
  Here $H$ is a minimal normal subgroup of $\Gamma$ and 
  an abelian $p$-group and $\pi_1 ^c(C)$ is the commutator of
  the algebraic fundamental group of a smooth affine curve $C$ over an
  algebraically closed field $k$ of characteristic $p$.
\end{thm}

\proof
Since $Z(\Gamma)$ the center of $\Gamma$ and $H$ are both normal 
subgroup of $\Gamma$, so is $Z(\Gamma)\cap H$. Since $H$ is minimal 
normal subgroup of $\Gamma$, $Z(\Gamma)\cap H$ is trivial or 
$H \subset Z(\Gamma)$. If $H \subset Z(\Gamma)$ then $\Gamma$ acts 
trivially on $H$, hence $\Gamma \cong G\times H$. By previous lemma there 
are $m$ surjections of $\pi_1 ^c(C)$
onto $H$. Borrowing the notation from Theorem 5.3, we conclude that
there are $m$ different field extensions of $K^{ab}$ contained in $K^{un}$
with Galois group $H$. $M$ being a finite field extension of $K^{ab}$ only 
finitely many of these $H$-extensions are not linearly disjoint with $M$ over $K^{ab}$. 
Hence there are $m$ $H$-extensions of $K^{ab}$ which are 
linearly disjoint with $M$ over $K^{ab}$ and compositum of each of these $H$-extensions
with $M$ lead to a $\Gamma$-extension of $K^{ab}$. Hence we have 
$m$ solution to the embedding problem. Now suppose $Z(\Gamma)\cap H$
is trivial, i.e. $\Gamma$ acts on $H$ nontrivially. By Theorem 5.1 there are 
$m$ proper solutions to the embedding problem
for $\pi_1(X)$ where $X$ is as in the proof of Theorem 5.3. So there are $m$
distinct smooth irreducible \'etale $H$-cover of $W_X$ which are $\Gamma$-covers 
of $X$.
For each such $H$-cover $Z$, we shall show that $\Gal(k(Z)K^{ab}/K^{ab})$ is isomorphic to
$\Gamma$. Suppose not, then $k(Z)$ is not linearly
disjoint with $M=k(W_X)K^{ab}$ over $k(W_X)$. So there exists a nontrivial field extension 
$L''/k(W_X)$ with $L'' = k(Z) \cap k(W_X)K^{ab}$. So $L''=Kk(W_X)$
for $K$ some finite field extension of $k(X)$ and $K\subset K^{ab}$. $K^{ab}/k(X)$ is a pro-abelian 
extension, so $K/k(X)$ is a Galois extension (with in fact abelian Galois group). Hence 
$L''/k(X)$ is a Galois extension. So we conclude that $\Gal(k(Z)/L'')$ is a normal subgroup
of $\Gamma=\Gal(k(Z)/k(X))$, but $\Gal(k(Z)/L'')\subset 
H=\Gal(k(Z)/k(W_X)$. $H$ being minimal normal subgroup of $\Gamma$ 
and $L''/k(W_X)$ being nontrivial extension, forces $L''=k(Z)$ and hence
$\Gal(K/k(X))=H$. But this contradicts the fact that $\Gamma$ acts 
on $H$ nontrivially. Now if $Z$ and $Z'$ are two distinct $H$-covers, replace 
the field $k(Z)$ by $k(Z)k(Z')$ and $\Gamma$ by $\Gal(k(Z)k(Z')/k(X))$ in 
the above argument to conclude that $k(Z)k(Z')$ is linearly disjoint with $M$ over $k(W_X)$.
Hence $k(Z)K^{ab}$ and $k(Z')K^{ab}$ are distinct field extensions of $M$.
\qed

Below we give an alternative approach to asserting existence of atleast one
proper solution of the embedding problem for $\pi_1 ^c(C)$ when
$H$ is any $p$-group. This is a cohomological approach.

\begin{thm}
  The following split embedding problem has a proper solution
  $$
  \xymatrix{
    &          &               &\pi_1 ^c(C) \ar[d] \ar @{-->} [dl]\\
    1\ar[r] & H \ar[r] & \Gamma \ar[r] & G \ar[r]\ar[d] & 1\\
    & & & 1 }
  $$
  Here $H$ is a $p$-group and $\pi_1 ^c(C)$ is the commutator of
  the algebraic fundamental group of a smooth affine curve $C$ over an
  algebraically closed field $k$ of characteristic $p$.
\end{thm}

\proof This theorem will follow trivially from Theorem 5.2, once we
prove that $H^1(\pi_1 ^c(C),P)$ is infinite for every elementary
abelian $p$-group $P$ with continuous $\pi_1^c(C)$-action, since we
have already seen that $\pi_1 ^c(C)$ is a projective profinite group, hence 
its cohomological dimension is less than or equal to $1$.  $H^1(\pi_1
^c(C),P)$ is infinite is shown in Proposition 5.7 below.  \qed

\begin{pro}
  Let $P$ be any nonzero finite elementary abelian $p$-group with a
  continuous action of $\pi_1 ^c(C)$, then the first group cohomology
  $H^1(\pi_1 ^c(C),P)$ is infinite.
\end{pro}

\proof Let $\Phi$ be the kernel of the action of $\pi_1 ^c(C)$ on $P$.
Then $\Phi$ is a normal subgroup of $\pi_1 ^c(C)$ of finite index. We
know $\pi_1 ^c(C)$ acts on the $K^{un}$ and has fixed field $K^{ab}$.
Let $M$ be the fixed field of $\Phi$, then $\Gal(M/K^{ab})=\pi_1
^c(C)/ \Phi$. Since $M$ is a finite extension of $K^{ab}$, there exists
$L$ a finite abelian extension of $k(C)$ and $L'$ a  finite extension
of $L$ such that $\Gal(L'/L)=\Gal(M/K^{ab})$ and $L'K^{ab}=M$. Let $X$
be the normalization of $C$ in $L$ and $Y$ be the normalization of $C$
in $L'$.  If we translate this Galois theory to Galois groups, we get
the following commutative diagram:
$$
\xymatrix {
  \Phi \ar@{^{(}->}[r] \ar@{^{(}->}[d] & \pi_(Y) \ar@{^{(}->}[d] \\
  \pi_1 ^c(C) \ar@{^{(}->}[r] \ar@{->>}[d]   & \pi_1(X) \ar@{->>}[d]\\
  \pi_1 ^c(C)/\Phi \ar[r]^{\sim} & \pi_1(X)/\pi_1(Y) }
$$
So we notice that $\pi_1(X)=\pi_1 ^c(C)\pi_1(Y)$. Now we define an
action of $\pi_1(X)$ on $P$ by defining it to be trivial on $\pi^1(Y)$
and to be the given action on $\pi_1 ^c(C)$.  This is well defined
because $\pi_1 ^c(C) \cap \pi_1(Y) = \Phi$, which is the kernel of the action
of $\pi_1 ^c(C)$ on
$P$.  Now consider the following short exact sequence of groups:
$$1 \rightarrow \pi_1 ^c(C) \rightarrow \pi_1(X) \rightarrow \Pi
\rightarrow 1$$
where $\Pi$ is simply the quotient $\pi_1(X)/\pi_1
^c(C)$. Applying Hoschild-Serre spectral sequence for group cohomology
[Wei 7.5.2] to this short exact sequence, we get the following long
exact sequence:
$$0 \rightarrow H^1(\Pi, H^0(\pi_1^c(C), P)) \rightarrow H^1(\pi_1^
c(X), P) \rightarrow H^0(\Pi, H^1(\pi_1 ^c(C),P))$$
$$\rightarrow H^2(\Pi,H^0(\pi_1^c(C),P))$$
If the action of $\pi_1
^c(C)$ on $P$ is such that it fixes only $0$, then $H^0(\pi_1 ^c(C),
P)=0$, hence the first and the fourth term in the above long exact
sequence is $0$. Also we know that the $H^1(\pi_1 ^c(X),P)$ is
infinite by [Ha4, Proposition 3.8].  So the $H^1(\pi_1 ^c(C),P)
\supset H^0(\Pi, H^1(\pi_1 ^c(C),P))$ is infinite.  So we may assume
$P^{\pi_1 ^c(C)}$ is nonzero.  In this case we have a short exact
sequence of $\pi_1^c(C)$-modules:
$$0 \rightarrow P^{\pi_1 ^c(C)} \rightarrow P \rightarrow P/P^{\pi_1
  ^c(C)} \rightarrow 0$$
Here $\pi_1 ^c(C)$ acts trivially on the
first term and fixes nothing in the third term, i.e., $H^0(\pi_1
^c(C), P/P^{\pi_1 ^c(C)})=0$, so we get the long exact sequence of
group cohomology which looks like:
$$... \rightarrow H^0(\pi_1 ^c(C), P/P^{\pi_1 ^c(C)})=0 \rightarrow
H^1(\pi_1 ^c(C), P^{\pi_1 ^c(C)}) \rightarrow H^1(\pi_1 ^c(C),P) ...$$
Since $\pi_1 ^c(C)$ acts trivially on $P^{\pi_1 ^c(C)}$, $H^1(\pi_1
^c(C), P^{\pi_1 ^c(C)})=\Hom(\pi_1 ^c(C), P^{\pi_1 ^c(C)})$.  And we
know that $\Hom(\pi_1 ^c(C), P^{\pi_1 ^c(C)})$ is infinite by previous
lemma.  So we conclude that $H^1(\pi_1 ^c(C),P)$ is infinite.  \qed
\\
{\bf Remark:} This alternative method could be possibly made strong enough to yield $m$ solutions
if one could prove a refined version of Theorem 5.2 relating the cardinality of $H^1(\pi,P)$
with the cardinality of distinct proper solutions for $p$-group embedding problem for $\pi$.
This looks plausible since the two method seems similar in spirit.

\section{Prime-to-$p$ embedding problem}

In this section, we prove certain results on formal patching and then
use them to solve the prime-to-$p$ embedding problems for $(*)$. We begin
with some patching results (6.1, 6.2, 6.3) which roughly mean
the following: suppose we have a proper $k[[t]]$-scheme $T$ whose
special fiber is a collection of smooth irreducible curves
intersecting at finitely many points. Finding a cover of $T$ is
equivalent to finding a cover of these irreducible curves away from
those finitely many intersection points and covers of formal
neighbourhood of the intersection points so that they agree in the
punctured formal neighbourhoods of the intersection points. In our
situation, the special fiber of $T$ is connected sum of $X$ and $N$
copies of $Y$.  Each copy of $Y$ intersects $X$ at a point $r_i$ of
$X$ and a point $s$ of $Y$ for $1\le i\le n$. Now suppose we have an
irreducible $G$-cover, $\Psi_X: W_X \rightarrow X$ \'etale at
$r_1,..,r_n$ and an irreducible $H$-cover, $\Psi_Y: W_Y \rightarrow Y$
\'etale at $s$ then we construct a $\Gamma$ covering of $T$ by
patching a $\Gamma$-cover of $X$, $\Ind^{\Gamma} _G W_X = (\Gamma
\times W_X)/\sim$, where $(\gamma,w) \sim(\gamma g^{-1}, gw)$ for
$\gamma\in\Gamma$, $g\in G$ and $w$ a point of $W_X$, and a
$\Gamma$-cover of $Y$, $\Ind^{\Gamma} _H W_Y$.  This is possible since
both these covers restrict to $\Gamma$-covers induced from trivial
cover in the formal punctured neighbourhood of the intersection
points, so we can pick trivial $\Gamma$-covers of the intersection
points which obviously will restrict to trivial $\Gamma$-cover
on the punctured neighbourhood. Now we proceed to show how all this
works. We start with the following patching result.

\begin{thm}({\bf [Ha3, Theorem 3.2.12]})
Let $(A,\cp)$ be a complete local ring and let $T$ be a proper $A$-scheme.
Let $\{\tau_1,..,\tau_N\}$ be a set of closed points of $T$ and 
$T^o=T\setminus\{\tau_1,..,\tau_N\}$. Let $\hat{T_i}=\spec(\hat{\cO}_{T,\tau_i})$, 
$\widetilde{T^o}$ be the $\cp$-adic completion of $T^o$ and $\cK_i$ be the 
$\cp$-adic completion of $\hat{T_i}\setminus\{\tau_i\}$. Then the base change functor
$$\cP(T)\rightarrow \cP(\widetilde{T^o})\times_{\cP(\cup_{i=1} ^N \cK_i}) 
\cP(\cup_{i=1} ^N \hat{T_i})$$ is an equivalence of categories. And this remains true
with $\cP$ replaced by $\cAP$, $\cSP$ or $G\cP$ for any finite group $G$.
\end{thm}

In fact [Ha3, Theorem 3.2.12] is even stronger and allows one to assert the 
equivalence of categories even if one replaces $T$, $T^o$, etc. with their 
pull back by a proper morphism. The proof uses Grothendieck's Existence Theorem and 
a result of Ferrand-Raynaud or rather its genralization by M. Artin 
([Ha3, Theorem 3.1.9]). 
The latter asserts, for a noetherian scheme $T$, the equivalence of categories between 
$\cP(T)$ and $\cP(T^o)\times_{\cP(W^o)} \cP(\hat{W})$ where $W$ is a finite set of 
closed points of $T$, $T^o=T\setminus W$, $\hat{W}$ is the completion of $T$ along $W$
and $W^o=\hat{W}\times_T T^o$.

Now we shall specialize to what we need. Let $k$ be a field.
Let $X$ and $Y$ be irreducible smooth projective $k$-curves with finite
$k$-morphisms $\Phi_X:X\rightarrow\PP^1 _x$ and
$\Phi_Y:Y\rightarrow\PP^1 _y$, where $\PP^1 _x$ and $\PP^1 _y$ are
projective lines with local coordinate $x$ and $y$ respectively. Also assume
that $\Phi_Y$ is totally ramified at $y=0$.  Let $R$ and $S$ be such
that $\spec(R)=X\setminus\Phi_X ^{-1}(\{x=\infty\})$ and
$\spec(S)=Y\setminus\Phi_Y ^{-1}(\{y=\infty\})$. So $k[x]\subset R$
and $k[y] \subset S$.  Let $A=(R\otimes_kS\otimes_kk[[t]])/(t-xy)$ and
$T^a=\spec(A)$. Let $T$ be the closure of $T^a$ in
$X\times_kY\times_k\spec(k[[t]])$. Let $L$ be an affine line
$\spec(k[z])$.  The $k$-algebra homomorphism $k[[t]][z]\rightarrow A$
given by $z \mapsto x+y$ induces a $k[[t]]$-morphism $\phi$ from $T^a$
to $L^*=L\times_kk[[t]]$. Let $\lambda \in L$ be the closed point
$z=0$. $L$ is contained in $L^*$ as the special fiber, so $\lambda$
when viewed as a closed point of $L^*$ corresponds to the maximal
ideal $(z,t)$ in $k[[t][z]$.  Let
$\phi^{-1}(\lambda)=\{\tau_1,..,\tau_N\} \subset T^a$. Note that the
special fiber of $T$ is a reducible curve consisting of $X$ and $N$
copies of $Y$, each copy of $Y$ intersecting $X$ at $\tau_1,..,\tau_N$
since locus of $t=0$ is same as $xy=0$ in $T$ and the locus of $t=0$
and $x+y=0$ is same as the locus of $x=0$ and $y=0$. Let $r_i$ denote
the point of $X$ corresponding to $\tau_i$, so $\Phi_X
^{-1}(x=0)=\{r_1,..,r_N\}$ and $s$ denote the point on each copy of
$Y$ corresponding to $\tau_i$, so $s$ is the unique point of $Y$ lying
above $y=0$.  Borrowing notation from the previous lemma, let $T^o=T
\setminus \{\tau_1,..\tau_N\}$ and $X^o=X\setminus\{r_1,..,r_N\}$.
Let $\hat{T_i}=\spec(\hat{\cO}_{T,\tau_i})$ and Let $T_X=T^o \setminus
\{x=0\}$ which is the same as the closure of $\spec(A[1/x])$ in
$X^o\times_kY\times_k\spec(k[[t]])$. Similarly, define $T_Y=T^o
\setminus \{y=0\}$.

Recall that $\hat{K}_{X, r_i}$ is the quotient field of $\hat{\cO}_{X,
  r_i}$. Define $\cK_{X,r_i}=
\spec(\hat{K}_{X,r_i}[[t]]\otimes_{k[y]}\cO_{Y,s})$ where we regard
$\hat{K}_{X,r_i}[[t]]$ as $k[y]$-module via the homomorphism which
sends $y$ to $t/x$.  Similarly, define $\cK_Y
^i=\spec(\hat{K}_{Y,s}[[t]]\otimes_{\spec(k[x])}\cO_{X,r_i})$, where
we regard $\hat{K}_{Y,s}[[t]]$ as $k[x]$-module via the homomorphism
which sends $x$ to $t/y$.  Let $x_i$ be a local coordinate of $X$ at
$r_i$ and $y_0$ be a local coordinate of $Y$ at $s$.

With these notations we shall deduce the following result from Theorem 6.1.
This result is analogous to [Ha2, Corollary 2.2].

\begin{lm}
  The base change functor $$\cP(T)\rightarrow \cP(\widetilde{T_X} \cup
  \widetilde{T_Y}) \times_{\cP(\D\cup_{i=1} ^N(\cK_{X,r_i} \cup \cK_Y
    ^i))} \cP(\cup_{i=1} ^N\hat{T_i})$$
  is an equivalence of
  categories. Moreover, same assertion holds if one replaces $\cP$ by
  $\cAP$, $\cSP$ and $G\cP$ for a finite group $G$.
\end{lm}

\proof 
First of all we observe that
the closed fiber of $T^o$, which is the subscheme defined by the ideal
$(t)$, is disconnected. Since closed fiber of $T_X \cup T_Y$ is the
closed fiber of $T^o$ and the closed fibers of $T_X$ and $T_Y$, as a 
subset of the closed fiber of $T^o$, are open and
disjoint.  So if we consider their $(t)$-adic completion we get
$\widetilde{T^o}=\widetilde{T_X} \cup \widetilde{T_Y}$. Similarly the 
punctured spectrum $\hat{T_i}\setminus \{\tau_i\}$ is the spectrum of the
ring $k[[x_i,y_0]][(x+y)^{-1}])$. Since the only prime ideals 
of $k[[x_i,y_0]][(x+y)^{-1}]$ containing $(t)$ are $(x_i)$ and
$(y_0)$, we may first localize $k[[x_i,y_0]][(x+y)^{-1}]$ with respect to
the complement of $(x_i)\cup(y_0)$ then take the $(t)$-adic completion.
Now using [Mat, 8.15], we get that the $(t)$-adic completion of 
$\hat{T_i}\setminus \{\tau_i\}$ is $\cK_{X,r_i}\cup\cK_Y ^i$.
\qed

Let $G$ and $H$ be subgroups of a finite group $\Gamma$, such that
$\Gamma=G\rtimes H$.

\begin{pro}
  Under the notation and assumption of previous lemmas, let
  $\Psi_X:W_X\rightarrow X$ be an irreducible normal $G$-cover \'etale
  over $r_1,..,r_N$.  and $\Psi_Y:W_Y\rightarrow Y$ be an irreducible
  normal $H$-cover \'etale over $s$. Let $W_{XT}$ be the normalization
  of an irreducible dominating component of $W_X\times_X T$ and
  similarly $W_{YT}$ be the normalization of an irreducible dominating
  component of $W_Y\times_Y T$.  Then there exists an irreducible
  normal $\Gamma$-cover $W\rightarrow
  T$ such that \\
  (1 ) $W\times_T \widetilde{T_X}= \Ind^{\Gamma} _G \widetilde{W_{XT}\times_T T_X}$\\
  (1') $W\times_T \widetilde{T_Y}= \Ind^{\Gamma} _H \widetilde{W_{YT}\times_T T_Y}$\\
  (2 ) $W\times_T \hat{T}_i$ is a $\Gamma$-cover of $\hat{T}_i$ induced from the trivial cover. \\
  (3 ) $W\times_T \cK_{X,r_i}$ is a $\Gamma$-cover of $\cK_{X,r_i}$ induced from the trivial cover. \\
  (4 ) $W\times_T \cK_Y ^i$ is a $\Gamma$-cover of $\cK_Y ^i$ induced from the trivial cover. \\
  (5 ) $W/H \cong W_{XT}$ as a cover of $T$.
\end{pro}

\proof Let $\widetilde{W_X}=\Ind^{\Gamma} _G \widetilde{W_{XT}\times_T
  T_X}$ and $\widetilde{W_Y}=\Ind^{\Gamma} _H
\widetilde{W_{YT}\times_T T_Y}$.  So $\widetilde{W_X}$ and
$\widetilde{W_Y}$ are $\Gamma$-covers of $\widetilde{T_X}$ and
$\widetilde{T_Y}$ respectively.  Hence their union, $\widetilde{W^o}$,
is an object of $\Gamma\cP(\widetilde{T_X} \cup \widetilde{T_Y})$.
Now for each $i$, $\widetilde{W_X}\times_{\widetilde{T_X}}\cK_{X,r_i}
=\Ind^{\Gamma} _G W_X\times_X \cK_{X,r_i}$. But $W_X\times_X
\cK_{X,r_i}$ is $\card(G)$ copies of $\cK_{X,r_i}$, since $W_X$ is
\'etale over $r_i$. And similarly,
$\widetilde{W_Y}\times_{\widetilde{T_Y}}\cK_Y ^i =\Ind^{\Gamma} _H
W_Y\times_Y \cK_Y ^i$ which is $\Gamma$ copies of $\cK_Y ^i$, since
$W_Y\rightarrow Y$ is \'etale over $s$. Hence $\widetilde{W^o}$
restricted to $\D\cup_{i=1} ^N(\cK_{X,r_i} \cup \cK_Y^i)$ is a
$\Gamma$-cover induced from the trivial cover.  Let $\hat{W}_i$ be a
$\Gamma$-cover of $\hat{T}_i$ induced from the trivial cover. Then
their union, $\hat{W}$, is an object in $\Gamma\cP(\cup_{i=1}
^N\hat{T_i})$ which when restricted to $\D\cup_{i=1} ^N(\cK_{X,r_i}
\cup \cK_Y^i)$ obviously is a $\Gamma$-cover induced from the trivial
cover. So after fixing an isomorphism between the two trivial
$\Gamma$-covers of $\D\cup_{i=1} ^N(\cK_{X,r_i} \cup \cK_Y^i)$, we can
apply the above patching lemma and obtain an object $W$ in
$\Gamma\cP(T)$ which induces the covers $\widetilde{W^o}$ and
$\hat{W}$ on $\widetilde{T^o}$ and $\cup_{i=1} ^N\hat{T_i}$
respectively. Hence we get conclusion $(1)$ to $(4)$ of the
proposition.  So it remains to prove $W$ is irreducible and normal and
conclusion $(5)$ holds.  For irreducibility of $W$ we note that $G$
and $H$ generate $\Gamma$.  Suppose $W$ is reducible. Consider
$\Gamma^o$, the stabilizer of the identity component of $W$. So $W$
has $\card(\Gamma/\Gamma^o)$ irreducible component. Since $G$ is the
stabilizer of the identity component of $\widetilde{W_X}$, and $H$ is
the stabilizer of the identity component of $\widetilde{W_Y}$, $G$ and
$H$ is contained in $\Gamma^o$. Hence $\Gamma^o=\Gamma$. Hence $W$ is
irreducible.  To show $W$ is normal it is enough to show that for each
closed point $\sigma$ of $T$, $W_{\sigma}=W\times_T
\spec(\hat{\cO}_{T,\sigma})$ is normal. If $\sigma = \tau_i$ for some
$i$ then $W_{\sigma}$ is isomorphic to copies of $\hat{T_i}$ and hence
is normal.  Otherwise $\sigma$ belongs to $T_X$ (or $T_Y$). So
$W_{\sigma}$ is isomorphic to $\Ind^{\Gamma} _G
\widetilde{W_{XT}\times_T T_X} \times_{T_X}
\spec(\hat{\cO}_{T_X,\sigma})$, which is a union of copies of
$\spec(\hat{\cO}_{W_{XT} \times_T T_X, \sigma ' })$, where $\sigma '$
are points of $W_{XT}\times_T T_X$ lying above $\sigma$.  But $W_{XT}
\times_T T_X$ is normal.  Similar argument holds in the case when
$\sigma \in T_Y$.  Next we shall show that $W/H$ and $W_{XT}$
restricts to same $G$-cover on the patches $\widetilde{T_X}$,
$\widetilde{T_Y}$ and $\hat{T_i}$ for all $i$. So conclusion $(5)$
will follow from the previous lemma's assertion about the equivalence
of categories.  Clearly, both $W/H$ and $W_{XT}$ restricts to trivial
$G$-cover of $\hat{T_i}$. Now, $W_{XT} \times_T \widetilde{T_X}=
W_{XT} \times_T \widetilde(T)\times_{\widetilde{T}}
\widetilde{T_X}=\widetilde{W_{XT}} \times_{\widetilde{T}}
\widetilde{T_X}$ and this is same as $\widetilde{W_{XT}\times_T T_X}$
since $T_X$ is an open subscheme of $T$.  On the other hand $W/H
\times_T \widetilde{T_X}=(W\times_T \widetilde{T_X})/H$.  But by
$(1)$, this is same as $\widetilde{W_{XT}\times_T T_X}$. Finally,
$W_{XT} \times_T \widetilde{T_Y}= \Ind^G _{\{e\}} \widetilde{T_Y}$,
since the image of $T_Y$ under the map $T\rightarrow X$, is the
generic point. So the $G$-cover $W_{XT} \rightarrow T$ is trivial over
the subscheme $T_Y$. And by $(1')$, $W/H\times_T \widetilde{T_Y}=
(\Ind^{\Gamma} _H \widetilde{W_{TY}})/H$ which is same as
$\Ind^{\Gamma/H} _{H/H} \widetilde{T_Y}$ since $W_Y/H=Y$.  But
$\Gamma/H$ is $G$.  \qed

\begin{lm}
  Let $T$, $X$ and $Y$ be as in previous lemma. Let $D$ be an irreducible
  smooth projective $k$-curve. Assume that $\Phi_X:X\rightarrow \PP^1 _x$
  factors through $D$, i.e., there exist $\Phi_X':X\rightarrow D$ and 
  $\Theta:D\rightarrow \PP^1 _x$ such that their composition is $\Phi_X$.
  Also assume $\Phi_Y$ and $\Phi_X'$ are abelian covers. For any 
  $k[[t]]$-scheme $V$, let $V^g$ denote the generic fiber.
  Then the morphism $\Phi_X'\times\Phi_Y\times Id_{\spec(k[[t]])}$ restricted to $T$
  from $T$ to its image in $D\times_k \PP^1 _y\times_k \spec(k[[t]])$
  induces an abelian cover 
  of projective $k((t))$-curves $T^g\rightarrow D\times_k \spec(k((t)))$.
\end{lm}

\proof We need to show that the function field, $k(T)$, of $T$ is an
abelian extension of the field $k(D)\otimes_k k((t))$.  Note that $k(T)$ is the
compositum of $L_1=k(X)\otimes_k k((t))$
and $L_2$. Here $L_2$ is the function field of a dominating irreducible component of 
$$(Y\times_k \spec(k((t)))) \times_{\PP^1 _y\times_k \spec(k((t)))} 
(D\times_k \spec(k((t))))$$ 
where the morphism $D\times_k \spec(k((t))) 
\rightarrow \PP^1 _y\times_k \spec(k((t)))$ is the composition of 
$D\times_k\spec(k((t))) \rightarrow \PP^1 _x\times_k\spec(k((t)))$ with morphism
$\PP^1 _x\times_k \spec(k((t))) \rightarrow \PP^1 _y\times_k \spec(k((t)))$
defined in local coordinates by sending $y$ to $t/x$.
Since $L_1$ is a base change of finite extension of $k(x)$ by
$k(D)\otimes_k k((t))$ and $L_2$ is a base change of finite extension 
of the subfield $k(t/x)$ by $k(D)\otimes_k k((t))$, we have 
$L_1\cap L_2=k(D)\otimes_k k((t))$. Hence $L_1$ and $L_2$
are linearly dijoint over $k(D)\otimes_k k((t))$.
Now $\Gal(L_1/k(D)\otimes_k k((t)))$
is isomorphic to $\Gal(k(X)/k(D))$ and $\Gal(L_2/k(D)\otimes_k k((t)))$ is
isomorphic to $\Gal(k(Y)/k(y))$. Hence these groups are abelian, since
the latter groups are so. Using the fact that the Galois group of
compositum of linearly disjoint Galois field extensions is the direct sum
of the two Galois groups, we get that $\Gal(k(T)/k(D)\otimes_k k((t)))$ is a
direct sum of abelian groups, hence is abelian.  \qed


We shall see a variation of the following result, which is a special
case of [Ha2, Proposition 2.6, Corollary 2.7].
\begin{pro}({\bf Harbater})
  Let $k$ be an algebraically closed field. Let $X_0 ^s$ be a
  smooth projective connected smooth $k$-curve.  Let
  $\zeta^1,..,\zeta^r \in X_0 ^s$. Let $X_0$ and $X_1$ be irreducible
  normal projective $k[[t]]$-curves. Suppose $X_1$ has geometrically
  smooth closed fiber.
  Let $\psi: X_1\rightarrow X_0$ be a $G$-cover
  with generic fiber $\psi^g: X_1 ^g\rightarrow X_0 ^g$. 
  Assume $X_0=X_0 ^s \times_k \spec(k[[t]])$ and $X_1 ^g$ 
  has genus at least $1$. Also assume
  $\psi^g$ is a smooth $G$-cover \'etale away from
  $\{\zeta_1,..,\zeta_r\}$ where $\zeta_j=\zeta^j\times_k k((t))\in
  X_0 ^g$ for $1\le j\le r$. Then there exist smooth connected
  $G$-cover $\psi^s:X_1^s\rightarrow X_0 ^s$ \'etale away from
  $\{\zeta^1,..,\zeta^r\}$.
\end{pro}

The proof of the following result is also similar to the one given in
[Ha2]. Though in [Ha2] the assumption that $X_1$ is a nonconstant family and hence the
assertion of existence of $m$ distinct solutions is not made, it is possible to do this 
as we shall see below.

\begin{pro} 
  Let $k$ be an algebraically closed field. Let $X_0,..,X_3$ be
  irreducible normal projective $k[[t]]$-curves and for $i>0$, $X_i$
  have generically smooth closed fibers. For $i=1,2$ and $3$, let
  $\psi_i:X_i\rightarrow X_{i-1}$ be proper surjective
  $k[[t]]$-morphisms and $\psi_i ^g:X_i ^g\rightarrow X_{i-1} ^g$ be
  the induced morphisms on the generic fibers. Assume $X_0 ^g=X_0 ^s
  \times_k k((t))$ for some smooth projective $k$-curves $X_0 ^s$ and 
  $X_1 ^g$ is of genus atleast $1$.
  Let $\zeta^1,..,\zeta^r \in X_0 ^s$ and $\zeta_j=\zeta^j\times_k
  k((t))\in X_0 ^g$ for $1\le j\le r$, so that $\psi_1 ^g\circ\psi_2
  ^g\circ\psi_3 ^g$ is \'etale away from $\{\zeta_1,..,\zeta_r\}$.
  Let $\psi_1$ be an $A$-cover, $\psi_2$ be a $G$-cover, $\psi_3$ be
  an $H$-cover and $\psi_2\circ\psi_3$ be a $\Gamma$-cover.  Then
  there exist $X_1 ^s$, $X_2 ^s$ and $X_3 ^s$ connected smooth
  projective $k$-curves and morphisms $\psi_i ^s:X_i ^s\rightarrow
  X_{i-1} ^s$ so that $\psi_1^s\circ\psi_2 ^s\circ\psi_3 ^s$ is
  \'etale away from $\{\zeta^1,..,\zeta^r\}$ and $\psi^1$ is an
  $A$-cover, $\psi^2$ is a $G$-cover, $\psi^3$ is an $H$-cover and
  $\psi^2\circ\psi^3$ is a $\Gamma$-cover. 
\end{pro}
\proof Since all the three groups are finite, the covers $\psi_i$ for
$i=1,..,3$ descends to a $B$-morphism, where $B \subset k[[t]]$ is a
regular finite type $k[t]$-algebra. That is, there exist connected
$B$-schemes $X_i ^B$ and morphism $\psi_i ^B:X_i ^B\rightarrow X_{i-1}
^B$ where $\psi_1 ^B$ is an $A$-cover, $\psi_2 ^B$ is a $G$-cover,
$\psi_3 ^B$ is an $H$-cover and $\psi_2 ^B\circ\psi_3 ^B$ is a
$\Gamma$-cover and $\psi_i ^B$ induces $\psi_i$. Moreover for
$E=\spec(B[t^{-1}])$, $X_i ^E=X_i ^B\times_B E$ are regular and $X_0
^E$ is isomorphic to $X_0 ^s\times_k E$. The induced morphism $\psi_i
^E$ are such that $\psi_1 ^E$ is an $A$-cover, $\psi_2 ^E$ is a
$G$-cover, $\psi_3 ^E$ is an $H$-cover, $\psi_2 ^E\circ\psi_3 ^E$ is a
$\Gamma$-cover and $\psi_1 ^E\circ\psi_2 ^E\circ\psi_3 ^E$ is ramified
only over $\{\zeta_E ^1,..\zeta_E ^r\}$.  To complete the proof, we
shall show that there exists a nonempty open subset $E'$ of $E$ so
that the fiber of $\psi_i ^E$ over each closed point of $E'$ is
irreducible and nonempty.  First we note that by [Ha2, Lemma 2.4(b)]
the closed fibers of $X_i \rightarrow \spec(k[[t]])$ are connected,
since by assumption the closed fibers are generically smooth. Hence
the fibers of $\psi_i ^B$ over $(t=0)$ are connected because $X_i ^B$
induces $X_i$.  Since $X_i$'s are normal, $X_i ^B$'s are unibranched
along the corresponding fibers over $(t=0)$. Hence by [Ha1,
Proposition 5], we have a nonempty open subset of $\spec(B)$, and
hence subset $E'$ of $E=\spec(B)\setminus(t=0)$, such that for all
closed points $e\in E'$ the fibers $X_i ^e$ of $X_i ^E \rightarrow E$ 
over $e$ are irreducible. Next, we shall show that there exist a nonempty open
subset $S$ of $E'$ such that the restriction morphism $X_i ^S \rightarrow S$
is smooth of relative dimension $1$. Since $k$ is algebraically closed $k(X_0 ^s)$ is 
seperably generated over $k$. Hence $k(X_0 ^{E'})$ is seperably generated over $k(E')$.
Moreover, since $\psi_i ^E$ are finite seperable morphisms (in fact their composition
is \'etale away from $\{\zeta_E ^1,..\zeta_E ^r\}$), we have 
$k(X_0 ^{E'})$ is seperably generated over $k(E')$. Since $X_i ^{E'}\rightarrow E'$
is a morphism integral schemes of relative dimesion $1$ and is generically seperable,
the relative sheaf of differentials is free of rank $1$ at the generic point ([Eis, Corollary 16.17a]). 
Hence there exist an open subset $S$ of $E'$ such that 
the morphism $X_i ^S\rightarrow S$ is smooth of relative dimension $1$. Moreover, 
the fiber $X_i ^s$  at each point $s\in S\subset E'$ is irreducible.
\qed

\begin{lm}
  There exists an abelian cover $Y\rightarrow \PP_y ^1$ ramified only
  at $y=0$, where it is totally ramified, with genus of $Y$
  arbitrarily large.
\end{lm}
\proof Let $Y'$ be the normal cover of $\PP_y ^1$ defined by the
equation $u^{p^n} - u - y^{p^{n}+1}=0$. To see this is an irreducible
polynomial in $k(y)[u]$, by Gauss lemma, it is enough to show it is
irreducible in $k[y,u]$. But in fact, it is irreducible in $k(u)[y]$
since $p^n+1^{th}$ root of $u^{p^n}-u$ does not belong to $k(u)$. Also
$Y'$ is \'etale everywhere except $y=\infty$ and since there is only
one point in $Y'$ lying above $\infty$ it is totally ramified there.
So by translation we can get $Y$, a cover of $\PP_y ^1$, which is
totally ramified at $y=0$ and \'etale elsewhere. Also the genus of
$Y$, by the Hurwitz formula, is given by the equation
$2g(Y)-2=(p^n+1)(g(\PP^1 _u)-2) + \deg(R)$ where $R$ is the
ramification divisor of the morphism $Y\rightarrow \PP^1 _u$.  We also
know that $\deg(R)=\sum_{P\in Y} e_P-1$ where $e_P$ is the
ramification index at the point $P \in Y$.  Now branch locus of $Y$ as
a cover of $\PP^1 _u$ is given by $u^{p^n}-u=0$ and $u=\infty$. For
each point $P\in Y$ lying above a branch point other than $\infty$,
$e_P=p^n+1$, so we get that $\deg(R) \ge p^np^n$. So we get the
following inequality for genus of $Y$.
$$2g(Y)-2 \ge -2(p^n+1)+p^{2n}$$
$$\Rightarrow g(Y)\ge p^n(p^n-2)/2$$
Clearly $g(Y)$ could be made
arbitrary large. Also note that $\Gal(k(Y)/k(y)) \cong \bigoplus_{i=1}
^n \ZZ/p\ZZ$. Hence $Y$ is an abelian cover of $\PP^1 _y$.  \qed

\begin{thm}
  The following split embedding problem has a proper solution
  $$
  \xymatrix{
    &          &               &\pi_1 ^c(C) \ar[d] \ar @{-->} [dl]\\
    1\ar[r] & H \ar[r] & \Gamma \ar[r] & G \ar[r]\ar[d] & 1\\
    & & & 1 }
  $$
  Here $H$ is prime to $p$-group minimal normal subgroup of $\Gamma$ and $\pi_1 ^c(C)$ is the
  commutator of the algebraic fundamental group of a smooth affine curve $C$
  over an algebraically closed field $k$ of characteristic
  $p$.
\end{thm}

\proof Let $K^{un}$ denote the compositum (in some fixed algebraic
closure of $k(C)$) of the function field of all Galois \'etale covers of $C$. 
And let $K^{ab}$ be the subfield of
$K^{un}$ obtained by considering only abelian covers with above
property. In these terms $\pi_1 ^c(C)$ is $\Gal(K^{un}/K^{ab})$.
So giving a surjection from $\pi_1 ^c(C)$ to $G$ is same as
giving a Galois extension $M\subset K^{un}$ of $K^{ab}$ with Galois
group $G$. Since $K^{ab}$ is an algebraic extension of $k(C)$ and $M$
is a finite extension of $K^{ab}$, we can find a finite abelian
extension $L \subset K^{ab}$ of $k(C)$ and $L'\subset K^{un}$ a
$G$-Galois extension of $L$ so that $M=K^{ab}L'$. Let $D$ be the smooth 
compactification of $C$, $X$ be the
normalization of $D$ in $L$ and $ \Phi_X':X\rightarrow D$ be the normalization
morphism. 
Then $X$ is an abelian cover of $D$ \'etale over
$C$ and the function field, $k(X)$, of $X$ is $L$. Let $W_X$ be
the normalization of $X$ in $L'$ and $\Psi_X$ corresponding
normalization morphism. Then $\Psi_X$ is \'etale away from points lying above
$D\setminus C$ and $k(W_X)=L'$. Since $k$ is algebraically closed, $k(C)/k$ has a 
seperating transedence basis. By a stronger version of Noether normalization (for instance, 
see [Eis, Corollary 16.18]), 
there exist a finite proper $k$-morphism from $C$ to $\Aff^1 _x$, where $x$ denotes the
local coordinate of the affine line, which is generically seperable. The branch 
locus of such a morphism is codimension $1$, hence this morphism is \'etale away from 
finitely many points. By translation we may assume none of these points map to $x=0$.
This morphism extends to a finite proper morphism $\Theta: D\rightarrow \PP^1 _x$.
Let $\Phi_X:X\rightarrow \PP^1 _x$ be the composition $\Theta\circ \Phi_X'$. 
Let $\{r_1,..,r_N\}=\Phi_X^{-1}(\{x=0\})$, then $\Phi_X$
is \'etale at $r_1,..,r_N$. Also note that $\Theta^{-1}(\{x=\infty\})=
D\setminus C$.
Let $l>0$ be any integer.
Let $\Phi_Y:Y\rightarrow \PP^1 _y$ be an abelian cover \'etale
everywhere except $y=0$, where it is totally ramified and genus of $Y$ 
is at least $2$ and more than the number of generators for $H^l$,
i.e., product of $H$ with itself $l$ times. 
Let $s$ be the point lying above $y=0$.
Existence of such a $Y$ is guaranteed by Lemma 6.7.
Since $H^l$ is prime-to-$p$, and $Y$ has high genus by [SGAI, XIII,
Corollary 2.12, page 392], there exists irreducible \'etale
$H^l$-cover $W_Y'$ of $Y$. By
taking appropriate quotient we get $l$ distinct
\'etale $H$-covers of $Y$. 
For $1\le i\le l$,
let $^i\Psi_{Y}:$$^iW_Y \rightarrow Y$ denote the
covering morphisms. Now we can apply proposition 6.3 for each $i$.  So we have an
irreducible normal $\Gamma$-cover $^iW\rightarrow T$ satisfying
conclusion $(1)$ to $(5)$ of the Proposition 6.3.  Also, by Lemma 6.4,
we know that the morphism $T\rightarrow B$, where $B$ is the locus of
$xy-t=0$ in $D \times_k\PP^1
_y\times_k\spec(k[[t]])$, induces an abelian cover of $D
\times_k \spec(k((t)))$. Let $V^g$ denote the
generic fiber of a $k[[t]]$-schemes $V$.  Since the branch locus of
$^iW ^g\rightarrow T^g$ is determined by the branch locus of $^iW \rightarrow T$ on
the patches. From $(1),(1')$ and $(2)$, we conclude that $^iW ^g\rightarrow T^g$ is
ramified only at points of $T^g$ lying above $x=\infty$ since 
$^iW_{YT}\rightarrow T_Y$ is \'etale everywhere and $W_{XT}\rightarrow T_X$
is \'etale away from the points which maps to $D \setminus C=\Theta^{-1}(\{x=\infty\})$ under the compostion
of the morphisms $W_{XT}\rightarrow T_X\rightarrow X\rightarrow D$.
Also $T^g
\rightarrow D \times_k \spec(k((t)))$ is ramified only at points above
$D\setminus C= \Theta^{-1}(\{x=\infty\})$, since $T \rightarrow B$ 
is ramified only at points above $x=\infty$ and
$y=0$ and on the generic fiber (i.e., $t \neq 0$) these two points get identified. 
So for each $1\le i\le l$, we get
the following diagram.
$$
\xymatrix{ ^iW^g \ar[d] \\
  ^iW^g/H \cong W_{XT} ^g \ar[d]^{G \text{ cover ramified only at pts lying above }x=\infty } \\
  T^g \ar[d]^{\text{abelian cover ramified only at }x=\infty }\\
  D\times \spec(k((t))) }
$$
$T$ clearly is not defined over $k$ and genus of $T^g$ is least genus of $Y$, hence at least 1. 
Now applying Proposition 6.6, 
we get the following diagram for each $i$, with same ramification properties as
above
$$
\xymatrix{ ^iW^s \ar[d] \\
  ^iW^s/H \cong W_{XT} ^s \ar[d]^{G\text{ cover ramified only at pts lying above }x=\infty }\\
  T^s \ar[d]^{\text{abelian cover ramified only at }x=\infty}\\
  D }
$$
where $-^s$, as in Proposition 6.6, denote the specialization to
the base field $k$. Note that $k(^iW^s)$ are linearly disjoint over $k(W_{XT} ^s)$ for $1\le i\le l$.
So to complete the proof, it is enough to show that for atleast one $i$, the
Galois group of $k(^iW^s)K^{ab}$ over $K^{ab}$ is $\Gamma$, where
$k(^iW^s)$ is the function field of $^iW^s$. 
Note that $k(W_X)\subset
k(W_{XT}^s) \subset k(W_X)k(T^s)$.  So
$k(W_{XT}^s)K^{ab}=k(W_X)K^{ab}$, since $k(T^s) \subset K^{ab}$.  By
assumption Galois group of $k(W_X)K^{ab}$ over $K^{ab}$ is $G$. So it
is enough to show that the Galois group of $k(^iW^s)K^{ab}$ over
$k(W_{XT}^s)K^{ab}$ is $H$ for some $i$.
Since $H$ is minimal normal subgroup of $\Gamma$, $H\cong \SSS\times\SSS\times..\times\SSS$,
for some simple group $\SSS$. If $\SSS$ is non abelian then $\SSS$ and hence
$H$ is perfect.
$\Gal(k(^iW^s)/k(W_{XT}^s))$ is perfect and
$k(W_{XT}^s)K^{ab}/k(W_{XT}^s)$ is a pro-abelian field extension, so they are
linearly disjoint. Hence $\Gal(k(^iW^s)K^{ab}/k(W_{XT}^sK^{ab}) \cong H$. 
If $\SSS$ is abelian then $\SSS \cong \ZZ/q\ZZ$ for some prime $q$
different from $p$. By Grothedieck's result on prime-to-$p$ part of the fundamental group 
(see Theorem 3.4), there are only finitely many nontrivial surjections from
$\pi_1(C)$ to $H$. These epimorphisms correspond to the $H$-covers $Z_j$ of $D$ which are \'etale 
over $C$.
Now note that we could have chosen $l$ to be any integer. So let $l$
be an integer greater that the number of such $H$-covers of $D$.
After base change, some of these $Z_j\times_D W_{XT}^s$ may still be $H$-covers of $W_{XT}^s$. 
We choose an 
$i$ such that $^iW^s$ is different from $Z_j\times_D W_{XT}^s$ for all $j$. For such an $i$,
$k(^iW^s)$ is linearly disjoint with $k(W_{XT}^s)K^{ab}$ over $k(W_{XT})$, since 
subfields of $k(W_{XT}^s)K^{ab}$ which are finite extensions of $k(W_{XT}^s)$ are in bijective 
correspondence with the covers of $W_{XT}^s$ obtained from base change of an \'etale cover of $C$.
So we found a proper solution to the embedding problem. 

\qed

\end{document}